\documentclass[reqno]{amsart}

\newcommand{\Rm}{\mathbb{R}}
\newcommand{\Cm}{\mathbb{C}}

\newcommand{\Sm}{\mathbb{S}}
\newcommand{\be}{\[}
\newcommand{\ee}{\]}
\newcommand{\ba}{\[\begin{aligned}}
\newcommand{\ea}{\end{aligned}\]}
\newcommand{\va}{\varphi}

\newcommand{\pp}{\partial}
\newcommand{\bv}[1]{\boldsymbol{\mathrm{#1}}}

\newcommand{\uv}{\boldsymbol{{\hat{\theta}}}}

\usepackage{amsmath}
\usepackage[foot]{amsaddr}
\usepackage{amsthm,amssymb,mathrsfs}
\usepackage{graphicx}

\theoremstyle{remark}

\title[]{The radiative transport equation with waiting time and its diffusion approximation with a time-fractional derivative}

\author[]{Manabu Machida}
\address{Department of Informatics, Faculty of Engineering, Kindai University, Higashi-Hiroshima 739-2116, Japan}
\email{machida@hiro.kindai.ac.jp}

\begin{document}

\begin{abstract}
Albeit the past intensive research, the governing equation of anomalous diffusion which is observed for the transport of particles underground is still an open problem. In this paper, as a governing equation, the advection-diffusion equation with a time-fractional derivative term is derived from the radiative transport equation with an integral term for waiting time.
\end{abstract}

\maketitle

\section{Introduction}

Anomalous diffusion has been observed in different phenomena \cite{Metzler-Klafter00}. One example is mass transport in underground water \cite{Adams-Gelhar92}. As another example, an equation with an additional time-fractional term \cite{Langlois-Farazmand-Haller15} is derived from the fluid dynamics in the Basset problem \cite{Basset88a,Basset88b,Basset88c,Basset10,Boussinesq85}. Aiming at reproducing experimental results for mass transport in porous media, various partial differential equations with fractional derivatives have been proposed. The fractional advection-dispersion equation with time- and spatial-fractional derivatives was proposed \cite{Benson-etal00}. Fractional derivative models were compared \cite{Sun-etal17,Wei-etal16} and parameters in fractional equations were estimated \cite{Chakraborty-etal09,Kelly-etal17}. Orders of fractional derivatives were taken to be variables \cite{Sun-etal09}. Furthermore, the distributed order time fractional diffusion equation was tested \cite{Liang-etal19}. Trying to reproduce such anomalous diffusion, a multi-scaling tempered fractional-derivative model was proposed \cite{Zhang-etal13}. The tempered anomalous diffusion model, which has a time-fractional derivative, was used to analyze column flow experiment \cite{Suzuki-etal16}.

By the comparison of the numerical solution of the radiative transport equation and the observed concentration of tracer particles in a column experiment, it was found that the governing equation at the mesoscopic scale (i.e., the propagation distance, which is the column height, is comparable to the transport mean free path) is the radiative transport equation \cite{Amagai-etal21,Amagai-etal20}.

In this paper, we show that an equation with a time-fractional derivative is asymptotically derived from the radiative transport equation with a waiting time. The derived equation governs transport in diffusive regime, in which the propagation distance is much larger than the transport mean free path.

The radiative transport equation, which is a linear Boltzmann equation, has been used in various fields. The diffusion approximation of the radiative transport equation has been intensively studied. In addition to the intuitive approach \cite{Ishimaru78,Duderstadt-Martin79}, The diffusion approximation was considered with the asymptotic expansion in the whole space \cite{Ryzhik-etal96} and in the presence of the boundary \cite{Larsen-Keller74}. Bensoussan, Lions, and Papanicolaou explored the diffusion approximation in terms of homogenization \cite{Bensoussan-Lions-Papanicolaou79}. The diffusion approximation in the presence of boundaries was treated for the propagation of elastic (seismic) waves \cite{Bal01}. Moreover, the Robin boundary condition, which is imposed on the diffusion equation, was considered with the diffusion approximation \cite{Bardos-Santos-Sentis84}.

The remainder of the paper is organized as follows. In Sec.~\ref{rte}, the radiative transport equation with waiting time is introduced. In Sec.~\ref{laplace}, the Laplace transform is taken to treat the integral term for traps. The diffusion approximation is considered in Sec.~\ref{da}. Numerical solutions to the two equations which are connected by the diffusion approximation are compared in Sec.~\ref{numerics}. Finally, concluding remarks are given in Sec.~\ref{concl}.

\section{The radiative transport equation with waiting time}
\label{rte}

Let $\bv{v}_0\in\Rm^3$ be the velocity of particles and $\bv{c}\in\Rm^3$ be advection. The direction of the propagation of a particle is given by $\uv=\bv{v}_0/|\bv{v}_0|$, which is a unit vector in $\Sm^2$. Let $\psi(\bv{r},\uv,\tau)$ be the angular density of particles at position $\bv{r}\in\Rm^3$ in direction $\uv\in\Sm^2$ at time $\tau>0$. The initial condition is imposed as
\be
\psi(\bv{r},\uv,0)=a(\bv{r},\uv),
\ee
where $a$ is the initial distribution of particles. Let $\Phi(\tau)$ be the survival probability and $\sigma_a,\sigma_s,\sigma_{\rm trap}$ be positive constants. We consider the following radiative transport equation in $\Rm^3$:
\begin{equation}
\begin{aligned}
&
P_0\psi(\bv{r},\uv,\tau)+\left(\sigma_a+\sigma_s+\sigma_{\rm trap}\right)\psi(\bv{r},\uv,\tau)=
\sigma_s\int_{\Sm^2}p(\uv\cdot\uv')\psi(\bv{r},\uv',\tau)\,d\uv'
\\
&+
\sigma_{\rm trap}\int_0^{\tau}w(\tau-\tau')\psi(\bv{r},\uv,\tau')\,d\tau'+
\sigma_{\rm trap}\Phi(\tau)a(\bv{r},\uv),
\end{aligned}
\label{rte1}
\end{equation}
where
\be
P_0\psi(\bv{r},\uv,\tau)=\left(\pp_{\tau}+\left(\bv{v}_0+\bv{c}\right)\cdot\nabla\right)\psi(\bv{r},\uv,\tau).
\ee
The scattering phase function satisfies $\int_{\Sm^2}p(\uv\cdot\uv')\,d\uv=1$. In this paper, $\sigma_a$, $\sigma_s$ for absorption and scattering are assumed to be constant. The waiting-time distribution $w$ is described below.

For column experiments, breakthrough curves were reproduced using (\ref{rte1}) without $\sigma_{\rm trap}$ \cite{Amagai-etal21,Amagai-etal20}. In the case of $\sigma_{\rm trap}=0$, we see that (\ref{rte1}) governs mass transport if the particle flow in complicated paths is regarded as the transport of particles which undergo scattering \cite{Williams92a,Williams92b,Williams93a,Williams93b,Amagai-etal20}. Although the column experiments were performed using uniformly random beads and sand, it can be expected that in field experiments tracer particles are trapped during the propagation or they propagate in an immobile zone \cite{Hatano-Hatano98,Schumer-etal03}. This motivates us to introduce terms with $\sigma_{\rm trap}$ in (\ref{rte1}). The coefficient $\sigma_{\rm trap}$ on the left-hand side of (\ref{rte1}) counts particles which are trapped and the integral term with $\sigma_{\rm trap}$ on the right-hand side means particles that were trapped at time $\tau'$ restart propagation at time $\tau$.

A trapped particle reenters the flow after some waiting time or escape time. The function $w$ in (\ref{rte1}) is the waiting time distribution. If each particle immediately restarts the propagation even when it is trapped, the waiting-time function is $w(\tau)=\delta(\tau)$, where $\delta(\cdot)$ is Dirac's delta function. In this case, terms with $\sigma_{\rm trap}$ on both sides of (\ref{rte1}) cancel. In general, $w(\tau)$ has a decaying behavior as $\tau$ grows. We have
\be
\Phi(\tau)=1-\int_0^{\tau}w(\tau')\,d\tau'.
\ee
The term $\sigma_{\rm trap}\Phi(\tau)a(\bv{r},\uv)$ on the right-hand side of (\ref{rte1}) is the source term due to particles which are in dead-end pores at $\tau=0$.

Recently, numerical simulation with dead-end pores was performed \cite{Bordoloi-etal22}. The probability density functions of escape time obtained from trajectories of particles implies the power-law behavior of $w(\tau)$. Different power-law behaviors of the waiting time function were numerically observed for different shapes of dead-end pores \cite{Hou-Jiang-Wu18}.

In light of these past researches, let us assume that $w$ asymptotically decays as $\tau^{-(1+\alpha)}$, where $\alpha>0$. Moreover we assume $0<\alpha<1$ for a long tail of $w$. Let $\gamma>0$ be a constant. Examples of such function $w$ include
\be
w(\tau)=\frac{\alpha/\gamma}{(1+\tau/\gamma)^{1+\alpha}},\quad
w(\tau)=\frac{\alpha(\gamma/\tau)^{\alpha}}{\tau\left(1+(\gamma/\tau)^{\alpha}\right)^2},\quad
w(\tau)=\alpha\gamma^{\alpha}\tau^{-(1+\alpha)}e^{-(\gamma/\tau)^{\alpha}}.
\ee
In all cases, $\int_0^{\infty}w(\tau)\,d\tau=1$. For large $\tau$, they behave as
\be
w(\tau)\sim\alpha\gamma^{\alpha}\tau^{-(1+\alpha)}.
\ee

Let us introduce a nondecreasing function $W(\tau)$ such that $dW/d\tau=w$ ($0<\tau<\infty$), $W\in[0,1)$, $W(0)=0$. Corresponding to the above examples of $w$, we have
\be
W(\tau)=1-\frac{1}{(1+\tau/\gamma)^{\alpha}},\quad
W(\tau)=\frac{1}{1+(\gamma/\tau)^{\alpha}},\quad
W(\tau)=e^{-(\gamma/\tau)^{\alpha}},\quad 0\le \tau<\infty.
\ee

\section{Laplace transform}
\label{laplace}

Let us investigate the solution to (\ref{rte1}) on a large time-scale. To this end, we scale $\tau$ as
\be
t=\epsilon\tau,
\ee
where $\epsilon>0$ is a small number. 

We write the solution of (\ref{rte1}) as \cite{Larsen-Keller74}
\be
\psi(\bv{r},\uv,\tau)=\psi^i(\bv{r},\uv,\tau)+\psi^{\rm IL}(\bv{r},\uv,\tau),
\ee
where $\psi^i$ is the interior part and $\psi^{\rm IL}$ is the initial layer part. In particular,
\be
a(\bv{r},\uv)=\psi^i(\bv{r},\uv,0)+\psi^{\rm IL}(\bv{r},\uv,0).
\ee
The interior part $\psi^i(\bv{r},\uv,\tau)=\psi_{\epsilon}^i(\bv{r},\uv,t)$ satisfies
\begin{equation}
\begin{aligned}
&
P_{0,\epsilon}\psi_{\epsilon}^i(\bv{r},\uv,t)+\left(\sigma_a+\sigma_s+\sigma_{\rm trap}\right)\psi_{\epsilon}^i(\bv{r},\uv,t)
=\sigma_s\int_{\Sm^2}p(\uv\cdot\uv')\psi_{\epsilon}^i(\bv{r},\uv',t)\,d\uv'
\\
&+
\frac{\sigma_{\rm trap}}{\epsilon}\int_0^{t/\epsilon}w\left(\frac{t-t'}{\epsilon}\right)\psi_{\epsilon}^i(\bv{r},\uv,t')\,dt'
+\sigma_{\rm trap}\Phi\left(\frac{t}{\epsilon}\right)a(\bv{r},\uv),
\end{aligned}
\label{rtei}
\end{equation}
where
\be
P_{0,\epsilon}\psi_{\epsilon}^i(\bv{r},\uv,t)=\left(\epsilon\pp_t+\left(\bv{v}_0+\bv{c}\right)\cdot\nabla\right)\psi_{\epsilon}^i(\bv{r},\uv,t).
\ee

We write the Laplace transform of a function $f_{\epsilon}(t)=f(\tau)$ as
\be
(\mathcal{L}f_{\epsilon})(s)=\int_0^{\infty}e^{-st}f_{\epsilon}(t)\,dt,
\quad s\in\Cm.
\ee

Let us take the Laplace transform for (\ref{rtei}):
\be
\begin{aligned}
&
(\mathcal{L}P_{0,\epsilon}\psi_{\epsilon}^i)(\bv{r},\uv,s)+
\left(\sigma_a+\sigma_s+\sigma_{\rm trap}\right)(\mathcal{L}\psi_{\epsilon}^i)(\bv{r},\uv,s)=
\sigma_s\int_{\Sm^2}p(\uv\cdot\uv')(\mathcal{L}\psi_{\epsilon}^i)(\bv{r},\uv',s)\,d\uv'
\\
&+
\sigma_{\rm trap}(\mathcal{L}w)(\epsilon s)(\mathcal{L}\psi_{\epsilon}^i)(\bv{r},\uv,s)+
\sigma_{\rm trap}\frac{1-(\mathcal{L}w)(\epsilon s)}{s}a(\bv{r},\uv),
\end{aligned}
\ee
where
\be
(\mathcal{L}P_{0,\epsilon}\psi_{\epsilon}^i)(\bv{r},\uv,s)=
\epsilon s(\mathcal{L}\psi_{\epsilon}^i)(\bv{r},\uv,s)-\epsilon\psi_{\epsilon}^i(\bv{r},\uv,0)+
\left(\bv{v}_0+\bv{c}\right)\cdot\nabla(\mathcal{L}\psi_{\epsilon}^i)(\bv{r},\uv,s).
\ee

We set
\be
d(\tau)=1-\left(\frac{\gamma}{\tau}\right)^{\alpha}.
\ee
We note that $W\sim d$ ($\tau\to\infty$) and for any $\lambda>0$, 
$d(\lambda \tau)/d(\tau)\to1$ as $\tau\to\infty$. We have for small $s$ \cite{Karamata31,Feller71},
\begin{equation}
(\mathcal{L}w)(s)\sim d\left(\frac{1}{s}\right)=
1-(\gamma s)^{\alpha},\quad s>0.
\label{w_laplace_asym}
\end{equation}
Since $\epsilon>0$ is small, we can replace $(\mathcal{L}w)(\epsilon s)$ with $1-(\gamma\epsilon s)^{\alpha}$.

\section{Diffusion approximation}
\label{da}

Assuming slow advection, large scattering, small absorption, small trapping rate, and large life time, we introduce $\epsilon$ as
\be
\bv{c}\to\epsilon\bv{c},\quad
\sigma_s\to\frac{\sigma_s}{\epsilon},\quad
\sigma_a\to\epsilon\sigma_a,\quad
\sigma_{\rm trap}\to\epsilon\sigma_{\rm trap},\quad
\gamma\to\frac{\gamma}{\epsilon}.
\ee
Then we have
\begin{equation}
\begin{aligned}
&
\epsilon s(\mathcal{L}\psi_{\epsilon}^i)(\bv{r},\uv,s)-\epsilon\psi_{\epsilon}^i(\bv{r},\uv,0)+
\left(\bv{v}_0+\epsilon\bv{c}\right)\cdot\nabla(\mathcal{L}\psi_{\epsilon}^i)(\bv{r},\uv,s)
\\
&+
\left(\epsilon\sigma_a+\frac{\sigma_s}{\epsilon}+\epsilon\sigma_{\rm trap}\right)(\mathcal{L}\psi_{\epsilon}^i)(\bv{r},\uv,s)=
\frac{\sigma_s}{\epsilon}\int_{\Sm^2}p(\uv\cdot\uv')(\mathcal{L}\psi_{\epsilon}^i)(\bv{r},\uv',s)\,d\uv'
\\
&+
\epsilon\sigma_{\rm trap}\left(1-(\gamma s)^{\alpha}\right)
(\mathcal{L}\psi_{\epsilon}^i)(\bv{r},\uv,s)+
\epsilon\sigma_{\rm trap}\gamma^{\alpha}s^{\alpha-1}a(\bv{r},\uv),
\end{aligned}
\label{laplaceteps}
\end{equation}
where we used (\ref{w_laplace_asym}). We will formally expand $\psi^i$ as
\be
\psi^i(\bv{r},\uv,\tau)=\psi_{\epsilon}^i(\bv{r},\uv,t)
=\sum_{n=0}^{\infty}\epsilon^n\psi_n(\bv{r},\uv,t).
\ee

Let us set $t_1=\tau/\epsilon$. Since $\psi^{\rm IL}$ is responsible for the initial rapid change of $\psi$, we require $\psi^{\rm IL}$ to satisfy
\be
\left(\pp_{t_1}+\epsilon\left(\bv{v}_0+\epsilon\bv{c}\right)\cdot\nabla\right)
\psi^{\rm IL}\left(\bv{r},\uv,\tau\right)=
\left(\mathcal{K}_{\epsilon}-I\right)\left(\epsilon^2\sigma_a+\sigma_s+\epsilon^2\sigma_{\rm trap}\right)\psi^{\rm IL}\left(\bv{r},\uv,\tau\right),
\ee
where $I$ is the identity, and for a function $f(\uv,\tau)$, operator $\mathcal{K}_{\epsilon}$ is defined as
\be
\begin{aligned}
\mathcal{K}_{\epsilon}f(\uv,\tau)
&=
\frac{\sigma_s/\epsilon}{\epsilon\sigma_a+\sigma_s/\epsilon+\epsilon\sigma_{\rm trap}}
\int_{\Sm^2}p(\uv\cdot\uv')f(\uv',\tau)\,d\uv'
\\
&+
\frac{\epsilon\sigma_{\rm trap}}{\epsilon\sigma_a+\sigma_s/\epsilon+\epsilon\sigma_{\rm trap}}
\int_0^{\tau}w\left(\tau-\tau'\right)f\left(\uv,\tau'\right)\,d\tau'.
\end{aligned}
\ee
The initial condition is given by
\be
\psi^{\rm IL}(\bv{r},\uv,0)=a(\bv{r},\uv)-\psi_0(\bv{r},\uv,0)-\sum_{n=1}^{\infty}\epsilon^n\psi_n(\bv{r},\uv,0).
\ee
In the limit of $\epsilon\to0$, we have
\be
\pp_{t_1}\psi^{\rm IL}(\bv{r},\uv,\tau)=
\left(\mathcal{K}_0-I\right)\sigma_s\psi^{\rm IL}(\bv{r},\uv,\tau),
\ee
where
\be
\mathcal{K}_0f(\uv,\tau)=\int_{\Sm^2}p(\uv\cdot\uv')f(\uv',\tau)\,d\uv'.
\ee
The initial condition is written as
\be
\psi^{\rm IL}(\bv{r},\uv,0)=a(\bv{r},\uv)-\psi_0(\bv{r},\uv,0).
\ee
The largest eigenvalue of $\mathcal{K}_0$ is $1$ \cite{Dautray-Lions,Larsen-Keller74}. Hence, $\psi^{\rm IL}$ does not vanish if the initial value is nonzero and has a component proportional to the corresponding eigenfunction, which is independent of $\uv$. This implies
\be
\psi_0(\bv{r},\uv,0)=\frac{1}{4\pi}\int_{\Sm^2}a(\bv{r},\uv)\,d\uv.
\ee
Thus, $\psi^{\rm IL}$ decays exponentially in time. Since $t_1=\epsilon^{-2}t$, the decay rate in $t$ is very large, being proportional to $\epsilon^{-2}$. The initial layer part $\psi^{\rm IL}$ is vanishingly small outside the initial layer of duration $O(\epsilon^2)$.

By considering terms of order $O(\epsilon^{-1})$ in (\ref{laplaceteps}), we have
\be
(\mathcal{L}\psi_0)(\bv{r},\uv,s)=\int_{\Sm^2}p(\uv\cdot\uv')(\mathcal{L}\psi_0)(\bv{r},\uv',s)\,d\uv'.
\ee
The above relation implies that $\psi_0$ is independent of $\uv$ and we can write
\be
\psi_0(\bv{r},\uv,t)=\frac{1}{4\pi}u(\bv{r},t).
\ee
We have
\be
u(\bv{r},0)=a_0(\bv{r}),
\ee
where
\be
a_0(\bv{r})=\int_{\Sm^2}a(\bv{r},\uv)\,d\uv.
\ee

Next we collect terms of order $O(1)$:
\be
\bv{v}_0\cdot\nabla(\mathcal{L}\psi_0)(\bv{r},\uv,s)+\sigma_s(\mathcal{L}\psi_1)(\bv{r},\uv,s)=
\sigma_s\int_{\Sm^2}p(\uv\cdot\uv')(\mathcal{L}\psi_1)(\bv{r},\uv',s)\,d\uv'.
\ee
Thus,
\be
\psi_1(\bv{r},\uv,t)=
-\frac{|\bv{v}_0|}{4\pi(1-g)\sigma_s}\uv\cdot\nabla u(\bv{r},t),
\ee
where
\be
g=\int_{\Sm^2}\uv\cdot\uv' p(\uv\cdot\uv')\,d\uv'.
\ee

Terms of order $O(\epsilon)$ satisfies the relation below:
\be
\begin{aligned}
&
s(\mathcal{L}\psi_0)(\bv{r},\uv,s)-\psi_0(\bv{r},\uv,0)+
\bv{v}_0\cdot\nabla(\mathcal{L}\psi_1)(\bv{r},\uv,s)+
\bv{c}\cdot\nabla(\mathcal{L}\psi_0)(\bv{r},\uv,s)
\\
&+
\sigma_a(\mathcal{L}\psi_0)(\bv{r},\uv,s)+
\sigma_s(\mathcal{L}\psi_2)(\bv{r},\uv,s)=
\sigma_s\int_{\Sm^2}p(\uv\cdot\uv')(\mathcal{L}\psi_2)(\bv{r},\uv',s)\,d\uv'
\\
&-
\sigma_{\rm trap}(\gamma s)^{\alpha}(\mathcal{L}\psi_0)(\bv{r},\uv,p)+
\sigma_{\rm trap}\gamma^{\alpha}s^{\alpha-1}a(\bv{r},\uv).
\end{aligned}
\ee
By integrating the above equation over $\uv$, we obtain
\begin{equation}
\begin{aligned}
&
s(\mathcal{L}u)(\bv{r},s)-u(\bv{r},0)-
\frac{|\bv{v}_0|^2}{3(1-g)\sigma_s}\Delta(\mathcal{L}u)(\bv{r},s)+
\bv{c}\cdot\nabla(\mathcal{L}u)(\bv{r},s)
+\sigma_a(\mathcal{L}u)(\bv{r},s)
\\
&=
-\sigma_{\rm trap}(\gamma s)^{\alpha}(\mathcal{L}u)(\bv{r},p)+
\sigma_{\rm trap}\gamma^{\alpha}s^{\alpha-1}a_0(\bv{r}).
\end{aligned}
\label{order2eq}
\end{equation}

Let us set
\be
\eta=\gamma^{\alpha}\sigma_{\rm trap},\quad
D_0=\frac{|\bv{v}_0|^2}{3(1-g)\sigma_s}.
\ee
We can rewrite (\ref{order2eq}) as
\begin{equation}
\begin{aligned}
&
s(\mathcal{L}u)(\bv{r},s)-u(\bv{r},0)-D_0\Delta(\mathcal{L}u)(\bv{r},s)+
\bv{c}\cdot\nabla(\mathcal{L}u)(\bv{r},s)
\\
&+
\sigma_a(\mathcal{L}u)(\bv{r},s)
+\eta\left(s^{\alpha}(\mathcal{L}u)(\bv{r},s)-s^{\alpha-1}a_0(\bv{r})\right)=0.
\end{aligned}
\end{equation}

For a function $f(t)$, the Caputo derivative $\pp_t^{\alpha}$ ($0<\alpha<1$) is given by \cite{Caputo67,Podlubny99}
\begin{equation}
\pp_t^{\alpha}f(t)=\frac{1}{\Gamma(1-\alpha)}
\int_0^t(t-t')^{-\alpha}\pp_{t'}f(t')\,dt',
\quad0<\alpha<1.
\label{alpha1}
\end{equation}
Here, $\Gamma(\cdot)$ is the gamma function. Then we have
\begin{equation}
\mathcal{L}(\pp_t^{\alpha}f)(s)=
s^{\alpha}(\mathcal{L}f)(s)-s^{\alpha-1}f(0).
\label{alpha2}
\end{equation}
When defined in (\ref{alpha1}), the function $f$ needs to satisfy certain conditions (e.g., $f$ has to be differentiable). See \cite{KRY} for the mathematical treatment of the Caputo derivative. We assume that $\pp_t^{\alpha}f$ exists.

Finally, we obtain
\begin{equation}
\left\{\begin{aligned}
&
\pp_tu(\bv{r},t)+\eta\pp_t^{\alpha}u(\bv{r},t)
-D_0\Delta u(\bv{r},t)+\bv{c}\cdot\nabla u(\bv{r},t)
+\sigma_a u(\bv{r},t)=0,\quad \bv{r}\in\Rm^3,\quad t>0,
\\
&
u(\bv{r},0)=a_0(\bv{r}),\quad \bv{r}\in\Rm^3.
\end{aligned}\right.
\label{tfde}
\end{equation}

\section{Numerical test}
\label{numerics}

Let us suppose that the advection $\bv{c}$ is negligibly small and the speed $|\bv{v}_0|$ is $1$. Moreover let us consider the isotropic scattering of $p=1/(4\pi)$. To see the mathematical structure, we consider the following simple one-dimensional case:
\begin{equation}
\left\{\begin{aligned}
&
\left(\pp_t+\mu\pp_x+\sigma_a+\sigma_s+\sigma_{\rm trap}\right)\psi(x,\mu,t)=
\frac{\sigma_s}{2}\int_{-1}^1\psi(x,\mu,t)\,d\mu
\\
&+
\sigma_{\rm trap}\int_0^tw(t-t')\psi(x,\mu,t')\,dt'+
\sigma_{\rm trap}\Phi(t)a(x,\mu),
\\
&\quad
-\infty<x<\infty,\quad-1\le\mu\le1,\quad t>0,
\\
&
\psi(x,\mu,0)=a(x,\mu),\quad-\infty<x<\infty,\quad-1\le\mu\le1,
\end{aligned}\right.
\label{1dRTE}
\end{equation}
where $\mu$ is the cosine of the polar angle. To be specific, we set
\be
a(x,\mu)=\delta(x),\quad
w(\tau)=\frac{\alpha/\gamma}{(1+\tau/\gamma)^{1+\alpha}},\quad
0<\alpha<1,\quad\gamma>0.
\ee
We note that in this case
\be
(\mathcal{L}w)(s)=\alpha e^{\gamma s}E_{1+\alpha}(\gamma s),
\ee
where $E_{1+\alpha}$ is the generalized exponential integral: 
$E_{1+\alpha}(z)=z^{\alpha}\int_z^{\infty}e^{-t}/t^{1+\alpha}\,dt$, $z\in\Cm\setminus(-\infty,0]$.

Let $(\mathcal{L}\psi)(x,\mu,s)$ be the Laplace transform of $\psi(x,\mu,t)$. Here, $\mathcal{L}\psi$ satisfies
\be
\left(\mu\pp_x+\sigma_t\right)(\mathcal{L}\psi)(x,\mu,s)=
\frac{\sigma_s}{2}\int_{-1}^1(\mathcal{L}\psi)(x,\mu,s)\,d\mu
+q(x,s)
\ee
for $-\infty<x<\infty$, $-1\le\mu\le1$, where
\ba
q(x,s)
&=
\left(\sigma_{\rm trap}(\mathcal{L}\Phi)(s)+1\right)\delta(x),
\\
\sigma_t(s)
&=
\sigma_a+\sigma_s+\left(\sigma_{\rm trap}(\mathcal{L}\Phi)(s)+1\right)s.
\ea
We note that the relation $(\mathcal{L}\Phi)(s)=[1-(\mathcal{L}w)(s)]/s$ was used for the last line of the above equation.

By employing the analytical discrete ordinates method (ADO) \cite{Barichello11,Barichello-Siewert99a,Barichello-Siewert99b,Barichello-Garcia-Siewert00,Barichello-Siewert02,Siewert-Wright99}, we approximate $\mathcal{L}\psi$ with $v$ which satisfies
\be
\left(\mu_i\pp_x+\sigma_t(s)\right)v(x,\mu_i,s)=
\frac{\sigma_s}{2}\sum_{j=1}^Nw_j\left[v(x,\mu_j,s)+v(x,-\mu_j,s)\right]+
q(x,s)
\ee
for $-\infty<x<\infty$. Here, $\mu_i,w_i$ ($i=1,2,\dots,N$) are abscissas and weights for the Gauss-Legendre quadrature, respectively. We label $\mu_i$ as $0<\mu_1<\cdots<\mu_N<1$. We can write
\be
v(x,\mu_i,s)=\left(\sigma_{\rm trap}(\mathcal{L}\Phi)(s)+1\right)\sum_{j=1}^N\left[G(x,\mu_i;0,\mu_j;s)+G(x,\mu_i;0,-\mu_j;s)\right].
\ee
Here, the fundamental solution $G(x,\mu_i;y,\mu_j;s)$ satisfies
\ba
\left(\mu_i\pp_x+\sigma_t(s)\right)G(x,\mu_i;y,\mu_j;s)
&=
\frac{\sigma_s}{2}\sum_{j=1}^Nw_j\left[G(x,\mu_i;y,\mu_j;s)+G(x,\mu_i;y,-\mu_j;s)\right]
\\
&+
\delta(x-y)\delta_{ij},
\ea
where $\delta_{ij}$ is the Kronecker delta.

To obtain $G(x,\mu_i;y,\mu_j;s)$, we consider the homogeneous equation:
\be
\left(\mu_i\pp_x+\sigma_t(s)\right)v_{\nu}(x,\mu_i,s)=
\frac{\sigma_s}{2}\sum_{j=1}^Nw_j\left[v_{\nu}(x,\mu_i,s)+v_{\nu}(x,\mu_i,s)\right]
\ee
with the separation constant $\nu$, which implicitly depends on $s$. We write $v_{\nu}$ as
\be
v_{\nu}(x,\mu_i,s)=\phi(\nu,\mu_i,s)e^{-x/\nu}.
\ee
We impose the normalization condition:
\be
\sum_{i=1}^Nw_i\left(\phi(\nu,\mu_i,s)+\phi(\nu,-\mu_i,s)\right)=1.
\ee
We note the relation $\phi(-\nu,\mu_i,s)=\phi(\nu,-\mu_i,s)$. Assuming that $\nu\neq\mu_i/\sigma_t(s)$, we obtain
\begin{equation}
\phi(\nu,\mu_i,s)=\frac{\sigma_s\nu}{2}\frac{1}{\sigma_t(s)\nu-\mu_i}.
\label{phifunc}
\end{equation}
The following orthogonality relation holds:
\be
\sum_{i=1}^Nw_i\mu_i\left[\phi(\nu,\mu_i,s)\phi(\nu',\mu_i,s)-
\phi(\nu,-\mu_i,s)\phi(\nu',-\mu_i,s)\right]=
\mathcal{N}(\nu,s)\delta_{\nu,\nu'},
\ee
where
\be
\mathcal{N}(\nu,s)=\sum_{i=1}^Nw_i\mu_i\left[\phi(\nu,\mu_i,s)^2-\phi(\nu,-\mu_i,s)^2\right].
\ee

To compute $\nu$, we introduce
\be
\Xi_{\pm}=\begin{pmatrix}
\pm\mu_1 &&& \\ & \pm\mu_2 && \\ && \ddots & \\ &&& \pm\mu_N
\end{pmatrix},
\quad
\Phi_{\pm}(\nu,s)=\begin{pmatrix}
\phi(\nu,\pm\mu_1,s) \\ \phi(\nu,\pm\mu_2,s) \\ \vdots \\ \phi(\nu,\pm\mu_N,s)
\end{pmatrix},
\quad
\{W\}_{ij}=w_j.
\ee
Then the homogeneous equation can be written as
\ba
\left(\sigma_t(s)I_N-\frac{1}{\nu}\Xi_+\right)\Phi_+(\nu,s)&=
\frac{\sigma_s}{2}W\left(\Phi_+(\nu,s)+\Phi_-(\nu,s)\right),
\\
\left(\sigma_t(s)I_N-\frac{1}{\nu}\Xi_-\right)\Phi_-(\nu,s)&=
\frac{\sigma_s}{2}W\left(\Phi_+(\nu,s)+\Phi_-(\nu,s)\right),
\ea
where $I_N$ is the identity matrix of dimension $N$. The above equation can be rewritten as
\be
\begin{pmatrix}\Xi_+^{-1} & \\ & \Xi_-^{-1}\end{pmatrix}
\left[\sigma_t(s)\begin{pmatrix}I_N&\\&I_N\end{pmatrix}-\frac{\sigma_s}{2}
\begin{pmatrix}W&W\\W&W\end{pmatrix}\right]
\begin{pmatrix}\Phi_+(\nu,s)\\ \Phi_-(\nu,s)\end{pmatrix}
=
\frac{1}{\nu}
\begin{pmatrix}\Phi_+(\nu,s)\\ \Phi_-(\nu,s)\end{pmatrix}.
\ee
The fact $\Phi_{\pm}(-\nu,s)=\Phi_{\mp}(\nu,s)$ implies that $-\nu$ is an eigenvalue if $\nu$ is an eigenvalue. Thus, $\nu=\pm\nu_n(s)$ ($n=1,2,\dots,N$) are computed as eigenvalues of the above eigenproblem ($\nu_n(s)>0$ assuming $\nu_n(s)$ are nonzero). The eigenvectors of the eigenproblem are not necessary since $\phi(\nu,\mu_i,s)$ are already explicitly given in (\ref{phifunc}). We note that, unlike the standard ADO for the equation with real coefficients, in the present case $\nu_n(s)$ are complex due to $\sigma_t(s)\in\Cm$. The fundamental solution is given by
\be
G(x,\mu_i;y,\mu_j,s)=\pm\sum_{\pm\Re\nu_n>0}\frac{w_j}{\mathcal{N}(\nu_n,s)}
\phi(\nu_n,\mu_i,s)\phi(\nu_n,\mu_j,s)e^{-(x-y)/\nu_n(s)},
\ee
where upper signs are chosen for $x>y$ and lower signs are used for $x<y$.

Let us consider
\be
U_{\rm RTE}(x,t;\alpha)=\mathcal{L}^{-1}\sum_{i=1}^Nw_i\left[v(x,\mu_i,s)+v(x,-\mu_i,s)\right],
\quad x>0,\quad t>0,
\ee
where $\mathcal{L}^{-1}$ denotes the inverse Laplace transform. We have
\be
(\mathcal{L}U_{\rm RTE})(x,s;\alpha)=
\left(\sigma_{\rm trap}(\mathcal{L}\Phi)(s)+1\right)
\sum_{\Re\nu_n>0}\frac{1}{\mathcal{N}(\nu_n,s)}e^{-x/\nu_n(s)}.
\ee

The particle density $U_{\rm RTE}(x,t;\alpha)$ is given by
\ba
U_{\rm RTE}(x,t;\alpha)
&=
\mathcal{L}^{-1}\left(\mathcal{L}U_{\rm RTE}\right)(x,t;\alpha)
\\
&=\frac{1}{2\pi i}\int_{\sigma-i\infty}^{\sigma+i\infty}e^{st}(\mathcal{L}U_{\rm RTE})(x,s;\alpha)\,ds,
\ea
where $\sigma\in\Rm$ is taken to be greater than the largest real part of any singularity. Assuming that $U_{\rm RTE}(x,t)$ is real, we have
\be
\Re(\mathcal{L}U_{\rm RTE})(x,\sigma+i\omega;\alpha)=
\int_0^{\infty}e^{-\sigma t}U_{\rm RTE}(x,t;\alpha)\cos(\omega t)\,dt,
\quad\omega\in\Rm.
\ee
By integrating both sides of the above equation over $\omega$ after multiplying $(2e^{\sigma t}/\pi)\cos\omega t$ ($t>0$), we obtain \cite{Davies-Martin79,Ganapol08}
\begin{equation}
U_{\rm RTE}(x,t;\alpha)=
\frac{2e^{\sigma t}}{\pi t}\int_0^{\infty}\Re(\mathcal{L}U_{\rm RTE})\left(x,\sigma+i\frac{\omega}{t};\alpha\right)\cos\omega\,d\omega.
\label{cosrepresentation}
\end{equation}

The inverse Laplace transform with (\ref{cosrepresentation}) can be performed numerically. We will use the double-exponential formula \cite{Ooura-Mori91,Ooura-Mori99}. Let us introduce
\be
\va(y)=\frac{y}{1-e^{-K\sinh{y}}},\quad K=6.
\ee
Here, $K=6$ was found by numerical experiments. Then the above integral can be expressed as
\be
U_{\rm RTE}(x,t;\alpha)=\frac{2e^{\sigma t}}{\pi t}\int_{-\infty}^{\infty}\cos\left(M\va(y)\right)\Re(\mathcal{L}U_{\rm RTE})\left(x,\sigma+\frac{iM}{t}\va(y);\alpha\right)M\frac{d\va(y)}{dy}\,dy,
\ee
where $\omega=M\va(y)$, $M>0$. We introduce $h=\pi/M$. The above integral can be numerically computed with the trapezoidal rule \cite{Sugihara97,Trefethen-Weideman14}. We have
\be
U_{\rm RTE}(x,t;\alpha)\approx\frac{2e^{\sigma t}}{t}\sum_{j=-J}^J\cos\left(M\va(y_j)\right)\Re(\mathcal{L}U_{\rm RTE})\left(x,\sigma+\frac{iM}{t}\va(y_j);\alpha\right)\left.\frac{d\va(y)}{dy}\right|_{y_j},
\ee
where $J$ is an integer and $y_j=jh+\pi/(2M)$. We note that $d\va/dy\to0$ double-exponentially as $y\to-\infty$. We see that $\va(y)\to y$ double exponentially as $y\to\infty$ and $\cos(M\va(jh+\pi/(2M)))\sim\cos(Mjh+\pi/2)=0$ for large $j$.

For the numerical calculation, we set $\alpha=1/2$, $\sigma=0.04$, $\sigma_a=10^{-9}\,{\rm min}^{-1}$, $\sigma_s=1\,{\rm min}^{-1}$, $N=30$, $k_{\rm max}=40$, $M=40$.

Next we consider the one-dimensional version of (\ref{tfde}):
\begin{equation}
\left\{\begin{aligned}
&
\left(\pp_t+\eta\pp_t^{\alpha}\right)U_{\rm DE}(x,t;\alpha)
-D_0\pp_x^2U_{\rm DE}(x,t;\alpha)+\sigma_a U_{\rm DE}(x,t;\alpha)=0,\quad 
x\in\Rm,\quad t>0,
\\
&
U_{\rm DE}(x,0;\alpha)=a_0(x),\quad x\in\Rm,
\end{aligned}\right.
\label{1dfDE}
\end{equation}
where the advection term was ignored ($\bv{c}=\bv{0}$) and
\be
a_0(x)=2\delta(x).
\ee
We note that $D_0=1/(3\sigma_s)$. The function
\be
V(k,s)=\int_{-\infty}^{\infty}e^{-ikx}(\mathcal{L}U_{\rm DE})(x,s;\alpha)\,dx,
\quad-\infty<k<\infty
\ee
satisfies
\be
V(k,s)=\frac{1+\eta s^{\alpha-1}}{s+\eta s^{\alpha}+D_0k^2+\sigma_a}
\int_{-\infty}^{\infty}e^{-ikx}a_0(x)\,dx
\ee

Since $\int_0^{\infty}e^{-by}\,dy=1/b$, we have
\be
V(k,s)=2(1+\eta s^{\alpha-1})\int_0^{\infty}e^{-(s+\eta s^{\alpha}+D_0k^2+\sigma_a)y}\,dy.
\ee
By a straightforward calculation we can show
\be
\mathcal{L}\Theta(t-y)g_{\alpha}\left(\frac{t-y}{(\eta y)^{1/\alpha}}\right)
(\eta y)^{-1/\alpha}=
e^{-sy}(\mathcal{L}g_{\alpha})\left((\eta y)^{1/\alpha}s\right)=
e^{-(s+\eta s^{\alpha})y},
\ee
where $\Theta(\cdot)$ is the Heaviside step function and $g_{\alpha}(\cdot)$ is such that $\mathcal{L}g_{\alpha}=e^{-s^{\alpha}}$. We have
\be
\mathcal{L}^{-1}e^{-s^{\alpha}}=
\sum_{j=0}^{\infty}\frac{(-1)^j}{j!}\mathcal{L}^{-1}s^{\alpha j}=
\sum_{j=0}^{\infty}\frac{(-1)^j}{j!}\frac{t^{-\alpha j-1}}{\Gamma(-\alpha j)}=
\frac{1}{t}F_{\alpha}\left(\frac{1}{t^{\alpha}}\right)=
\frac{\alpha}{t^{1+\alpha}}M_{\alpha}\left(\frac{1}{t^{\alpha}}\right),
\ee
where using the Mainardi-Wright function $W_{\lambda,\mu}(z)$,
\be
F_{\alpha}(z)=W_{-\alpha,0}(-z),\quad
M_{\alpha}(z)=W_{-\alpha,1-\alpha}(-z),\quad0<\alpha<1.
\ee
Hence \cite{Mainardi-Mura-Pagnini10},
\be
g_{\alpha}(t)=\frac{1}{t}F_{\alpha}\left(\frac{1}{t^{\alpha}}\right)=
\frac{\alpha}{t^{1+\alpha}}M_{\alpha}\left(\frac{1}{t^{\alpha}}\right).
\ee
In particular, we have
\be
M_{1/2}(z)=\frac{1}{\sqrt{\pi}}e^{-z^2/4}.
\ee
We note that \cite{Schumer-etal03}
\begin{equation}
\lim_{\eta\to0}g_{\alpha}\left(\frac{t-y}{(\eta y)^{1/\alpha}}\right)
(\eta y)^{-1/\alpha}=
\lim_{\eta\to0}\frac{\alpha\eta y}{(t-y)^{1+\alpha}}M_{\alpha}\left(\frac{\eta y}{(t-y)^{\alpha}}\right)
=\delta(t-y).
\label{ga_delta}
\end{equation}

Thus,
\ba
&
(\mathcal{L}U_{\rm DE})(x,s;\alpha)
\\
&=
\frac{1+\eta s^{\alpha-1}}{\pi}\int_0^{\infty}
\left(\mathcal{L}\Theta(t-y)g_{\alpha}\left(\frac{t-y}{(\eta y)^{1/\alpha}}\right)(\eta y)^{-1/\alpha}\right)
\int_{-\infty}^{\infty}e^{ikx}e^{-(D_0k^2+\sigma_a)y}\,dkdy
\\
&=
\frac{1}{\pi}\mathcal{L}\left(\delta(t)+\frac{\eta t^{-\alpha}}{\Gamma(1-\alpha)}\right)(s)
\\
&\times
\int_0^{\infty}\left(\mathcal{L}\Theta(t-y)g_{\alpha}\left(\frac{t-y}{(\eta y)^{1/\alpha}}\right)(\eta y)^{-1/\alpha}\right)
\sqrt{\frac{\pi}{yD_0}}e^{-\frac{x^2}{4yD_0}}e^{-\sigma_ay}\,dy,
\ea
where we used under the condition $\Re s>0$,
\be
1+\eta s^{\alpha-1}=(\mathcal{L}\delta(t))(s)+\frac{\eta}{\Gamma(1-\alpha)}(\mathcal{L}t^{-\alpha})(s).
\ee
By the inverse Laplace transform we obtain
\ba
U_{\rm DE}(x,t;\alpha)
&=
\frac{1}{\eta^{1/\alpha}\sqrt{\pi D_0}}
\int_0^t\left(\delta(t-t')+\frac{\eta(t-t')^{-\alpha}}{\Gamma(1-\alpha)}\right)
\\
&\times
\left(\int_0^{t'}g_{\alpha}\left(\frac{t'-y}{(\eta y)^{1/\alpha}}\right)y^{-\frac{1}{\alpha}-\frac{1}{2}}e^{-\frac{x^2}{4yD_0}}e^{-\sigma_ay}\,dy\right)\,dt'
\\
&=
\frac{\alpha\eta}{\sqrt{\pi D_0}}\int_0^t\left(\delta(t-t')+\frac{\eta(t-t')^{-\alpha}}{\Gamma(1-\alpha)}\right)
\\
&\times
\left(\int_0^{t'}\frac{\sqrt{y}}{(t'-y)^{1+\alpha}}M_{\alpha}\left(\frac{\eta y}{(t'-y)^{\alpha}}\right)e^{-\frac{x^2}{4yD_0}}e^{-\sigma_ay}\,dy\right)\,dt'.
\ea
In the limit of $\sigma_{\rm trap}\to0$, we recover the normal diffusion:
\begin{equation}
\lim_{\eta\to0}U_{\rm DE}(x,t;\alpha)=
\frac{1}{\sqrt{\pi D_0}}t^{-\frac{1}{2}}e^{-\frac{x^2}{4tD_0}}e^{-\sigma_at},
\label{normaldiff}
\end{equation}
where the relation (\ref{ga_delta}) was used.

In the case of $\alpha=1/2$, we obtain
\ba
U_{\rm DE}\left(x,t;\frac{1}{2}\right)
&=
\frac{\eta}{2\pi\sqrt{D_0}}\int_0^t\left(\delta(t-t')+\frac{\eta}{\sqrt{\pi}(t-t')^{1/2}}\right)
\\
&\times
\left(\int_0^{t'}\frac{\sqrt{y}}{(t'-y)^{3/2}}e^{-\frac{\eta^2y^2}{4(t'-y)}}e^{-\frac{x^2}{4D_0y}}e^{-\sigma_ay}\,dy\right)\,dt',
\ea
where we used $\Gamma(1/2)=\sqrt{\pi}$. For a numerical purpose, we change variables as $y=tt_1$ and then change $t_1=1-\tau_1^2$ in the first term. Similarly we change variables as $t'=tt_1$, $y=tt_1t_2$ and then change $t_2$ as $t_2=1-\tau_2^2$ in the second term. We obtain
\ba
U_{\rm DE}\left(x,t;\frac{1}{2}\right)
&=
\frac{\eta}{\pi\sqrt{D_0}}\int_0^1\frac{\sqrt{1-\tau_1^2}}{\tau_1^2}
e^{-\frac{\eta^2t\left(1-\tau_1^2\right)^2}{4\tau_1^2}}e^{-\frac{x^2}{4D_0t\left(1-\tau_1^2\right)}}e^{-\sigma_at(1-\tau_1^2)}\,d\tau_1
\\
&+
\frac{\eta^2\sqrt{t}}{\pi\sqrt{\pi D_0}}\int_0^1\frac{1}{\sqrt{1-t_1}}
\\
&\times
\left(\int_0^1\frac{\sqrt{1-\tau_2^2}}{\tau_2^2}
e^{-\frac{\eta^2tt_1\left(1-\tau_2^2\right)^2}{4\tau_2^2}}
e^{-\frac{x^2}{4D_0tt_1\left(1-\tau_2^2\right)}}
e^{-\sigma_att_1(1-\tau_2^2)}\,d\tau_2\right)\,dt_1.
\ea

\begin{figure}[ht]
\centering
\includegraphics[width=1.0\textwidth]{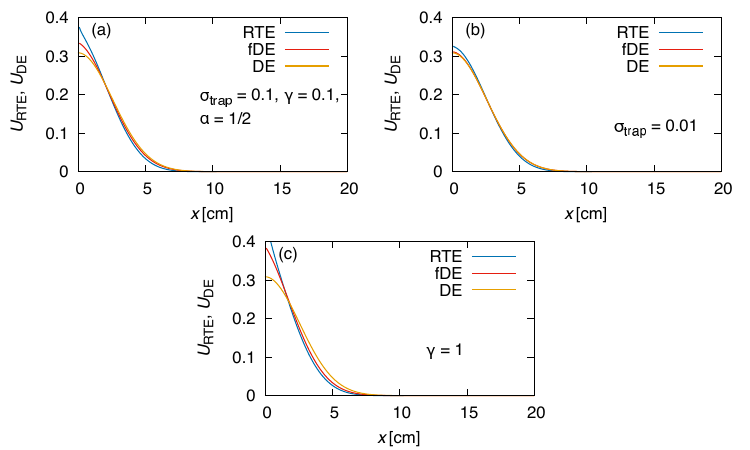}
\caption{
Comparison of the solution (blue line) to the transport equation (\ref{1dRTE}), the solution (red line) to the fractional diffusion equation (\ref{1dfDE}), and the solution (ocher line) to the normal diffusion (\ref{normaldiff}) at $t=10\,{\rm min}$. In panel (a), the parameters are set to $\sigma_{\rm trap}=0.1\,{\rm min}^{-1}$, $\gamma=0.1\,{\rm min}$, $\alpha=1/2$. The same parameters are used in other panels except that (b) $\sigma_{\rm trap}=0.01\,{\rm min}^{-1}$ and (c) $\gamma=1\,{\rm min}$, respectively.
}
\label{fig1}
\end{figure}

\begin{figure}[ht]
\centering
\includegraphics[width=1.0\textwidth]{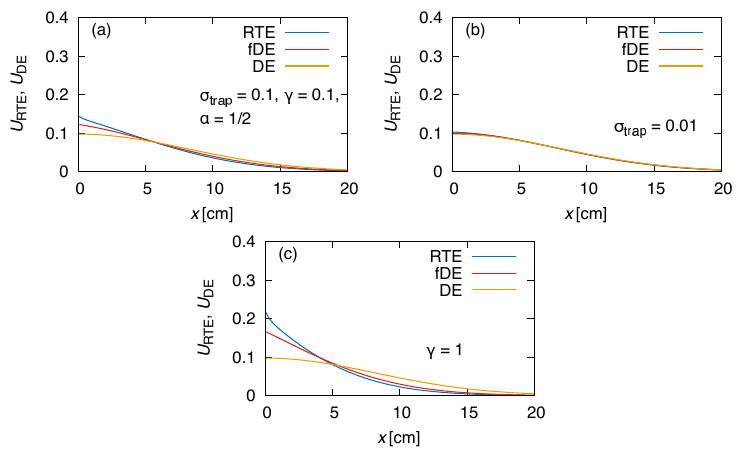}
\caption{
Same as Fig.~\ref{fig1} but $t=100\,{\rm min}$.
}
\label{fig2}
\end{figure}

Figures \ref{fig1} and \ref{fig2} show the solution $U_{\rm RTE}$ to the transport equation (\ref{1dRTE}) and the solution $U_{\rm DE}$ to the fractional diffusion equation (\ref{1dfDE}) as functions of $x\,[{\rm cm}]$ together with the normal diffusion (\ref{normaldiff}). We see that $U_{\rm RTE}$ and $U_{\rm DE}$ match as $x$ increases. The solutions at $t=10\,{\rm min}$ are shown in Fig.~\ref{fig1} and solutions at $t=100\,{\rm min}$ are shown in Fig.~\ref{fig2}. In panel (a), the parameters are set to $\sigma_{\rm trap}=0.1\,{\rm min}^{-1}$, $\gamma=0.1\,{\rm min}$, $\alpha=1/2$. In panel (b), the same parameters are used but $\sigma_{\rm trap}=0.01\,{\rm min}^{-1}$ to show that for a small $\sigma_{\rm trap}$, the propagation enters the diffusion regime for a shorter distance. In panel (c), $\gamma=1\,{\rm min}$ while other parameters are the same as those given in panel (a). When $\gamma$ is large, $U_{\rm RTE}$ and $U_{\rm DE}$ differ more for small $x$. Mathematica (Wolfram Research) was used to produce Figs.~\ref{fig1}, \ref{fig2}.

\section{Concluding remarks}
\label{concl}

Indeed, the derived fractional equation is the time-fractional advection-dispersion equation \cite{Schumer-etal03,Suzuki-etal16}. As a related model, the tempered anomalous diffusion model is known \cite{Meerschaert-etal08}. The fractional equation in this paper is derived from the radiative transport equation in a straightforward manner.

The waiting time distribution is also used for the continuous-time random walk. However, the physical meaning of the solution to (\ref{rte1}) is clearer than random walkers. That is, $\psi$ is the angular density of tracer particles while the relation between tracer molecules and Monte Carlo particles is vague.

If $\sigma_{\rm trap}$ is tiny (there is almost no dead-end pore), the fractional-derivative term in (\ref{tfde}) is negligible. In this case, the usual advection-diffusion equation is recovered. Finally, we note that the diffusion approximation can similarly be done in the case of the exponential decay, $w(\tau)=\frac{1}{\tau}e^{-\tau/\gamma}$, but the resulting diffusion equation only has a first-order time-derivative term.

\section*{Acknowledgements}

This work was supported by JST, PRESTO Grant Number JPMJPR2027.


\end{document}